\documentclass{article}

\usepackage{amsmath, amsfonts, amssymb, amsthm}
\usepackage{tikz}
\usepackage{url}
\usetikzlibrary{patterns}

\theoremstyle{plain}
\newtheorem{theorem}{Theorem}[section]
\newtheorem{proposition}[theorem]{Proposition}
\newtheorem{lemma}[theorem]{Lemma}
\newtheorem{definition}[theorem]{Definition}

\title{Every square can be tiled with T-tetrominos and no more than 5 monominos}
\author{Jack Grahl}

\begin{document}
\maketitle

\begin{abstract}
If $n$ is a multiple of 4, then a square of side $n$ can be tiled with T-tetrominos, using a well-known construction.
If $n$ is even but not a multiple of four, then there exists an equally well-known construction for tiling a square of side $n$ with T-tetrominos and exactly 4 monominos.
On the other hand, it was shown by Walkup in \cite{walkup} that it is not possible to tile the square using only T-tetrominos.
Now consider the remaining cases, where $n$ is odd.
It was shown by Zhan in \cite{zhan} that it is not possible to tile such a square using only one monomino.
Hochberg showed in \cite{hochberg} that no more than 9 monominos are ever needed.
We give a construction for all odd $n$ which uses exactly 5 monominos, thereby resolving this question.
\end{abstract}

\section{Introduction}
The sequence \cite{oeis} gives the maximal number of T-tetrominos which can be used to tile the $n \times n$ square with t-tetrominos and monominos. Theorem \ref{main} shows that this sequence is trivially given by $\frac{n^2}{4}, \frac{(n^2-1)}{4}-1, \frac{n^2}{4}-1, \frac{(n^2-1)}{4}-1$, depending on the value of $n$ modulo 4.
\section{Tiling every square}
\begin{theorem}\label{main}
Every square can be tiled with T-tetrominos and at most 5 monominos.
\end{theorem}
This theorem follows immediately from propositions \ref{four}, \ref{even} and \ref{odd}.

\begin{proposition}\label{four}
Every square of side $n = 4m$ can be tiled with T-tetrominos.
\end{proposition}
\begin{proposition}\label{even}
Every square of side $n = 4m + 2$ can be tiled with T-tetrominos and 4 monominos, and 4 monominos are always needed.
\end{proposition}
For $n = 2$ this is the same as pointing out that a single T-tetromino will not fit in the $2x2$ square.

\begin{figure}
\begin{tikzpicture}[scale=0.5]

\draw [rounded corners, ultra thick ] (0.020000, 0.500000) -- (0.020000, 0.980000) -- (0.980000, 0.980000) -- (0.980000, 0.020000) -- (0.020000, 0.020000) -- (0.020000, 0.500000);

\draw [rounded corners, ultra thick ] (0.020000, 1.500000) -- (0.020000, 1.980000) -- (0.980000, 1.980000) -- (0.980000, 1.020000) -- (0.020000, 1.020000) -- (0.020000, 1.500000);

\draw [rounded corners, ultra thick ] (1.020000, 0.500000) -- (1.020000, 0.980000) -- (1.980000, 0.980000) -- (1.980000, 0.020000) -- (1.020000, 0.020000) -- (1.020000, 0.500000);

\draw [rounded corners, ultra thick ] (5.020000, 5.500000) -- (5.020000, 5.980000) -- (5.980000, 5.980000) -- (5.980000, 5.020000) -- (5.020000, 5.020000) -- (5.020000, 5.500000);

\draw [rounded corners, ultra thick ] (2.500000, 1.980000) -- (3.980000, 1.980000) -- (3.980000, 1.020000) -- (2.980000, 1.020000) -- (2.980000, 0.020000) -- (2.020000, 0.020000) -- (2.020000, 1.020000) -- (1.020000, 1.020000) -- (1.020000, 1.980000) -- (2.500000, 1.980000);

\draw [rounded corners, ultra thick ] (4.500000, 0.020000) -- (3.020000, 0.020000) -- (3.020000, 0.980000) -- (4.020000, 0.980000) -- (4.020000, 1.980000) -- (4.980000, 1.980000) -- (4.980000, 0.980000) -- (5.980000, 0.980000) -- (5.980000, 0.020000) -- (4.500000, 0.020000);

\draw [rounded corners, ultra thick ] (5.980000, 2.500000) -- (5.980000, 1.020000) -- (5.020000, 1.020000) -- (5.020000, 2.020000) -- (4.020000, 2.020000) -- (4.020000, 2.980000) -- (5.020000, 2.980000) -- (5.020000, 3.980000) -- (5.980000, 3.980000) -- (5.980000, 2.500000);

\draw [rounded corners, ultra thick ] (4.020000, 4.500000) -- (4.020000, 5.980000) -- (4.980000, 5.980000) -- (4.980000, 4.980000) -- (5.980000, 4.980000) -- (5.980000, 4.020000) -- (4.980000, 4.020000) -- (4.980000, 3.020000) -- (4.020000, 3.020000) -- (4.020000, 4.500000);

\draw [rounded corners, ultra thick , pattern = north east lines] (0.020000, 4.000000) -- (0.020000, 5.980000) -- (3.980000, 5.980000) -- (3.980000, 2.020000) -- (0.020000, 2.020000) -- (0.020000, 4.000000);

\end{tikzpicture}
\caption{Extending the $4 \times 4$ tiling to $6 \times 6$, adding 4 monominos and 4 T-tetrominos.}
\label{six}
\end{figure}
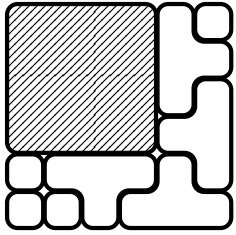

For $n = 4m + 2$, where $m$ is a positive integer, we can extend the tiling of the $4m$-square without monominos to a tiling of the $4m+2$-square, adding only 4 monominos. The tiling of the the L-shaped strip which extends the $4 \times 4$ square to a $6 \times 6$ square is given in figure \ref{six}. We can increase the length of the arms of the strip, by replacing the two T-tetrominos with a longer sequence taken from the `frieze', or tiling of a strip of width 2.

\begin{proposition}\label{odd}
Every square of side $n = 2m + 1$ can be tiled with T-tetrominos and 5 monominos, and 5 monominos are always needed (except for $n = 1$).
\end{proposition}
Zhan's (\cite{zhan}) Theorem 2 states that it is not possible to tile any rectangle with T-tetrominos and only one monomino. It must therefore be the case that at least 5 are needed. We show that exactly 5 are sufficient.

\begin{definition}
Call $A_n$ the set of lattice squares given by the square of side $n$, with the lattice squares at $(0, 0), (0,1), (1, 0)$ and $(0, n-1)$ removed.
This shape has area $n^2 - 4 = 4(m^2 + m - 1) + 1$.
\end{definition}

\begin{lemma}
For all $m \in \mathbb{N}$, $A_{2m+1}$ can be tiled with $m^2 + m - 1$ T-tetrominos and one monomino.
\end{lemma}

\begin{figure}
\begin{tikzpicture}[scale=0.5]

\draw [rounded corners, ultra thick ] (3.000000, 0.020000) -- (3.980000, 0.020000) -- (3.980000, 0.980000) -- (4.980000, 0.980000) -- (4.980000, 4.980000) -- (0.020000, 4.980000) -- (0.020000, 2.020000) -- (1.020000, 2.020000) -- (1.020000, 1.020000) -- (2.020000, 1.020000) -- (2.020000, 0.020000) -- (3.000000, 0.020000);

\draw [rounded corners, ultra thick ] (9.000000, 0.020000) -- (11.980000, 0.020000) -- (11.980000, 0.980000) -- (12.980000, 0.980000) -- (12.980000, 6.980000) -- (6.020000, 6.980000) -- (6.020000, 2.020000) -- (7.020000, 2.020000) -- (7.020000, 1.020000) -- (8.020000, 1.020000) -- (8.020000, 0.020000) -- (9.000000, 0.020000);

\draw [rounded corners, ultra thick ] (17.000000, 0.020000) -- (21.980000, 0.020000) -- (21.980000, 0.980000) -- (22.980000, 0.980000) -- (22.980000, 8.980000) -- (14.020000, 8.980000) -- (14.020000, 2.020000) -- (15.020000, 2.020000) -- (15.020000, 1.020000) -- (16.020000, 1.020000) -- (16.020000, 0.020000) -- (17.000000, 0.020000);

\end{tikzpicture}
\caption{$A_5$, the $5 \times 5$ square with four lattice squares removed, $A_7$ and $A_9$.}
\label{cropped}
\end{figure}
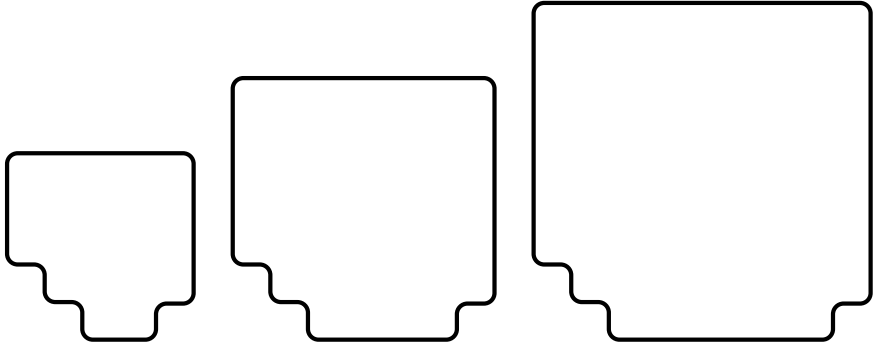

\begin{figure}
\begin{tikzpicture}[scale=0.5]

\draw [rounded corners, ultra thick ] (3.000000, 0.020000) -- (3.980000, 0.020000) -- (3.980000, 0.980000) -- (4.980000, 0.980000) -- (4.980000, 4.980000) -- (0.020000, 4.980000) -- (0.020000, 2.020000) -- (1.020000, 2.020000) -- (1.020000, 1.020000) -- (2.020000, 1.020000) -- (2.020000, 0.020000) -- (3.000000, 0.020000);

\draw [rounded corners, ultra thick ] (0.020000, 3.500000) -- (0.020000, 4.980000) -- (0.980000, 4.980000) -- (0.980000, 3.980000) -- (1.980000, 3.980000) -- (1.980000, 3.020000) -- (0.980000, 3.020000) -- (0.980000, 2.020000) -- (0.020000, 2.020000) -- (0.020000, 3.500000);

\draw [rounded corners, ultra thick ] (2.500000, 4.980000) -- (3.980000, 4.980000) -- (3.980000, 4.020000) -- (2.980000, 4.020000) -- (2.980000, 3.020000) -- (2.020000, 3.020000) -- (2.020000, 4.020000) -- (1.020000, 4.020000) -- (1.020000, 4.980000) -- (2.500000, 4.980000);

\draw [rounded corners, ultra thick ] (2.980000, 1.500000) -- (2.980000, 0.020000) -- (2.020000, 0.020000) -- (2.020000, 1.020000) -- (1.020000, 1.020000) -- (1.020000, 1.980000) -- (2.020000, 1.980000) -- (2.020000, 2.980000) -- (2.980000, 2.980000) -- (2.980000, 1.500000);

\draw [rounded corners, ultra thick ] (3.020000, 1.500000) -- (3.020000, 2.980000) -- (3.980000, 2.980000) -- (3.980000, 1.980000) -- (4.980000, 1.980000) -- (4.980000, 1.020000) -- (3.980000, 1.020000) -- (3.980000, 0.020000) -- (3.020000, 0.020000) -- (3.020000, 1.500000);

\draw [rounded corners, ultra thick ] (4.980000, 3.500000) -- (4.980000, 2.020000) -- (4.020000, 2.020000) -- (4.020000, 3.020000) -- (3.020000, 3.020000) -- (3.020000, 3.980000) -- (4.020000, 3.980000) -- (4.020000, 4.980000) -- (4.980000, 4.980000) -- (4.980000, 3.500000);

\draw [rounded corners, ultra thick , pattern = north east lines] (1.020000, 2.500000) -- (1.020000, 2.980000) -- (1.980000, 2.980000) -- (1.980000, 2.020000) -- (1.020000, 2.020000) -- (1.020000, 2.500000);

\end{tikzpicture}
\caption{Tiling of $A_5$ with a single monomino.}
\label{five}
\end{figure}
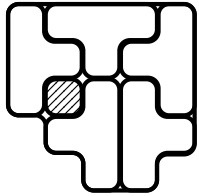

{\sc Proof.}
The proof is by induction on $n$. In figure \ref{five} we show how $A_5$ can be tiled by 5 tetrominos and a single monomino. (It is trivial to tile $A_3$ with a single tetromino and a single monomino, but it is slightly clearer to start the induction with $n=5$.) If $A_n$ can be tiled with one monomino, then so can $A_{n+1}$. There are two constructions for the cases $n=4k+1$ and $n=4k+3$.

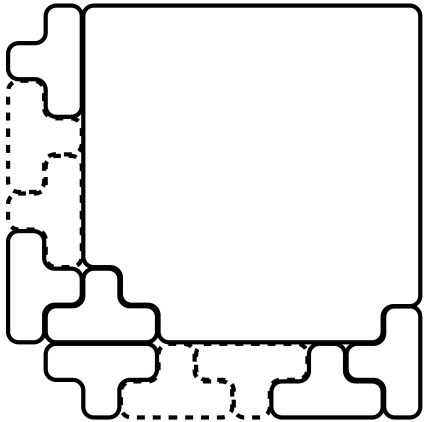
\begin{figure}
\begin{tikzpicture}[scale=0.5]

\draw [rounded corners, ultra thick ] (5.000000, 2.020000) -- (9.980000, 2.020000) -- (9.980000, 2.980000) -- (10.980000, 2.980000) -- (10.980000, 10.980000) -- (2.020000, 10.980000) -- (2.020000, 4.020000) -- (3.020000, 4.020000) -- (3.020000, 3.020000) -- (4.020000, 3.020000) -- (4.020000, 2.020000) -- (5.000000, 2.020000);

\draw [rounded corners, ultra thick ] (1.980000, 9.500000) -- (1.980000, 8.020000) -- (1.020000, 8.020000) -- (1.020000, 9.020000) -- (0.020000, 9.020000) -- (0.020000, 9.980000) -- (1.020000, 9.980000) -- (1.020000, 10.980000) -- (1.980000, 10.980000) -- (1.980000, 9.500000);

\draw [rounded corners, ultra thick ] (0.020000, 3.500000) -- (0.020000, 4.980000) -- (0.980000, 4.980000) -- (0.980000, 3.980000) -- (1.980000, 3.980000) -- (1.980000, 3.020000) -- (0.980000, 3.020000) -- (0.980000, 2.020000) -- (0.020000, 2.020000) -- (0.020000, 3.500000);

\draw [rounded corners, ultra thick ] (2.500000, 2.020000) -- (1.020000, 2.020000) -- (1.020000, 2.980000) -- (2.020000, 2.980000) -- (2.020000, 3.980000) -- (2.980000, 3.980000) -- (2.980000, 2.980000) -- (3.980000, 2.980000) -- (3.980000, 2.020000) -- (2.500000, 2.020000);

\draw [rounded corners, ultra thick ] (2.500000, 1.980000) -- (3.980000, 1.980000) -- (3.980000, 1.020000) -- (2.980000, 1.020000) -- (2.980000, 0.020000) -- (2.020000, 0.020000) -- (2.020000, 1.020000) -- (1.020000, 1.020000) -- (1.020000, 1.980000) -- (2.500000, 1.980000);

\draw [rounded corners, ultra thick ] (8.500000, 0.020000) -- (7.020000, 0.020000) -- (7.020000, 0.980000) -- (8.020000, 0.980000) -- (8.020000, 1.980000) -- (8.980000, 1.980000) -- (8.980000, 0.980000) -- (9.980000, 0.980000) -- (9.980000, 0.020000) -- (8.500000, 0.020000);

\draw [rounded corners, ultra thick ] (10.980000, 1.500000) -- (10.980000, 0.020000) -- (10.020000, 0.020000) -- (10.020000, 1.020000) -- (9.020000, 1.020000) -- (9.020000, 1.980000) -- (10.020000, 1.980000) -- (10.020000, 2.980000) -- (10.980000, 2.980000) -- (10.980000, 1.500000);

\draw [rounded corners, ultra thick , dashed] (0.020000, 7.500000) -- (0.020000, 8.980000) -- (0.980000, 8.980000) -- (0.980000, 7.980000) -- (1.980000, 7.980000) -- (1.980000, 7.020000) -- (0.980000, 7.020000) -- (0.980000, 6.020000) -- (0.020000, 6.020000) -- (0.020000, 7.500000);

\draw [rounded corners, ultra thick , dashed] (1.980000, 5.500000) -- (1.980000, 4.020000) -- (1.020000, 4.020000) -- (1.020000, 5.020000) -- (0.020000, 5.020000) -- (0.020000, 5.980000) -- (1.020000, 5.980000) -- (1.020000, 6.980000) -- (1.980000, 6.980000) -- (1.980000, 5.500000);

\draw [rounded corners, ultra thick , dashed] (4.500000, 0.020000) -- (3.020000, 0.020000) -- (3.020000, 0.980000) -- (4.020000, 0.980000) -- (4.020000, 1.980000) -- (4.980000, 1.980000) -- (4.980000, 0.980000) -- (5.980000, 0.980000) -- (5.980000, 0.020000) -- (4.500000, 0.020000);

\draw [rounded corners, ultra thick , dashed] (6.500000, 1.980000) -- (7.980000, 1.980000) -- (7.980000, 1.020000) -- (6.980000, 1.020000) -- (6.980000, 0.020000) -- (6.020000, 0.020000) -- (6.020000, 1.020000) -- (5.020000, 1.020000) -- (5.020000, 1.980000) -- (6.500000, 1.980000);

\end{tikzpicture}
\caption{A tiling of $A_{4k+1}$ can be extended to a tiling of a reflected copy of $A_{4k+3}$.}
\label{ones}
\end{figure}

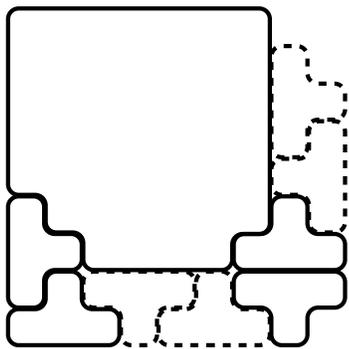
\begin{figure}
\begin{tikzpicture}[scale=0.5]

\draw [rounded corners, ultra thick ] (3.000000, 2.020000) -- (5.980000, 2.020000) -- (5.980000, 2.980000) -- (6.980000, 2.980000) -- (6.980000, 8.980000) -- (0.020000, 8.980000) -- (0.020000, 4.020000) -- (1.020000, 4.020000) -- (1.020000, 3.020000) -- (2.020000, 3.020000) -- (2.020000, 2.020000) -- (3.000000, 2.020000);

\draw [rounded corners, ultra thick ] (0.020000, 2.500000) -- (0.020000, 3.980000) -- (0.980000, 3.980000) -- (0.980000, 2.980000) -- (1.980000, 2.980000) -- (1.980000, 2.020000) -- (0.980000, 2.020000) -- (0.980000, 1.020000) -- (0.020000, 1.020000) -- (0.020000, 2.500000);

\draw [rounded corners, ultra thick ] (1.500000, 0.020000) -- (0.020000, 0.020000) -- (0.020000, 0.980000) -- (1.020000, 0.980000) -- (1.020000, 1.980000) -- (1.980000, 1.980000) -- (1.980000, 0.980000) -- (2.980000, 0.980000) -- (2.980000, 0.020000) -- (1.500000, 0.020000);

\draw [rounded corners, ultra thick ] (7.500000, 1.980000) -- (8.980000, 1.980000) -- (8.980000, 1.020000) -- (7.980000, 1.020000) -- (7.980000, 0.020000) -- (7.020000, 0.020000) -- (7.020000, 1.020000) -- (6.020000, 1.020000) -- (6.020000, 1.980000) -- (7.500000, 1.980000);

\draw [rounded corners, ultra thick ] (7.500000, 2.020000) -- (6.020000, 2.020000) -- (6.020000, 2.980000) -- (7.020000, 2.980000) -- (7.020000, 3.980000) -- (7.980000, 3.980000) -- (7.980000, 2.980000) -- (8.980000, 2.980000) -- (8.980000, 2.020000) -- (7.500000, 2.020000);

\draw [rounded corners, ultra thick , dashed] (3.500000, 1.980000) -- (4.980000, 1.980000) -- (4.980000, 1.020000) -- (3.980000, 1.020000) -- (3.980000, 0.020000) -- (3.020000, 0.020000) -- (3.020000, 1.020000) -- (2.020000, 1.020000) -- (2.020000, 1.980000) -- (3.500000, 1.980000);

\draw [rounded corners, ultra thick , dashed] (5.500000, 0.020000) -- (4.020000, 0.020000) -- (4.020000, 0.980000) -- (5.020000, 0.980000) -- (5.020000, 1.980000) -- (5.980000, 1.980000) -- (5.980000, 0.980000) -- (6.980000, 0.980000) -- (6.980000, 0.020000) -- (5.500000, 0.020000);

\draw [rounded corners, ultra thick , dashed] (8.980000, 4.500000) -- (8.980000, 3.020000) -- (8.020000, 3.020000) -- (8.020000, 4.020000) -- (7.020000, 4.020000) -- (7.020000, 4.980000) -- (8.020000, 4.980000) -- (8.020000, 5.980000) -- (8.980000, 5.980000) -- (8.980000, 4.500000);

\draw [rounded corners, ultra thick , dashed] (7.020000, 6.500000) -- (7.020000, 7.980000) -- (7.980000, 7.980000) -- (7.980000, 6.980000) -- (8.980000, 6.980000) -- (8.980000, 6.020000) -- (7.980000, 6.020000) -- (7.980000, 5.020000) -- (7.020000, 5.020000) -- (7.020000, 6.500000);

\end{tikzpicture}
\caption{A tiling of $A_{4k+3}$ can be extended to a tiling of a reflected copy of $A_{4(k+1)+1}$.}
\label{threes}
\end{figure}

\bibliography{polyomino}{}

\begin{thebibliography}{1}

\bibitem{oeis}
Jack Grahl.
\newblock Sequence {A}256535 of the {O}nline {E}ncyclopedia of {I}nteger
  {S}equences.
\newblock \url{http://oeis.org/A256535}, 2015.

\bibitem{hochberg}
Robert Hochberg.
\newblock The gap number of the {T}-tetromino.
\newblock 2014.

\bibitem{walkup}
D.~W. Walkup.
\newblock Covering a rectangle with {T}-tetrominos.
\newblock {\em The American Mathematical Monthly}, 72(9), November 1965.

\bibitem{zhan}
Shuxin Zhan.
\newblock Tiling a deficient rectangle with {T}-tetrominos.
\newblock 2012.

\end{thebibliography}
\bibliographystyle{plain}
\end{document}